\documentclass[journal,letterpaper,twocolumn,10pt,twosided]{IEEEtran}

\usepackage{amsmath}
\usepackage{amsmath,amssymb}
\usepackage{graphicx}
 \usepackage[utf8]{inputenc}   
\usepackage{color}
\usepackage{multirow}
\usepackage{hhline}
\usepackage{rotating,mathdots}
\usepackage{cite}

\usepackage{bbm}  
\usepackage{bm}  
\usepackage{booktabs} 
\usepackage{array} 
\usepackage{paralist} 
\usepackage{subfig} 
\usepackage{accents}
\usepackage[percent]{overpic}

\newcommand{\ubar}[1]{\underaccent{\bar}{#1}}

\newtheorem{proposition}{Proposition}

\newtheorem{consequence}{Consequence}

\makeatletter
\def\endexempel{\hfill\qed\@endtheorem}
\makeatother

\graphicspath{{Figures/}}
\DeclareGraphicsExtensions{.pdf,.jpg}

\newcommand*{\NN}{\mathbb{N}}
\newcommand*{\RR}{\mathbb{R}}
\newcommand*{\ZZ}{\mathbb{Z}}
\newcommand*{\EE}[1]{\mathbb{E}\!\left[#1\right]}

\newcommand*{\leadfini}{\breve\ell^{(p)}}

\newcommand*{\sfcoef}[1][q]{S_{\coe{}}(#1,j)}

\newcommand*{\sflead}[1][q]{S_{\ell^{(p)}}(#1,j)}
\newcommand*{\sfleadfini}[1][q]{S_{\leadfini}(#1,j)}
\newcommand*{\sfleadfinico}[1][q]{\hat S_{\leadfini}(#1,j)}

\newcommand*{\cumjlead}[2][j]{C_{\ell^{(p)}}(#2, #1)}
\newcommand*{\cumjleadfini}[2][j]{C_{\leadfini}(#2, #1)}
\newcommand*{\cumjleadfinico}[2][j]{\hat C_{\leadfini}(#2, #1)}

\newcommand*{\zetaq}{\zeta_{\ell^{(p)}}}
\newcommand*{\legspec}{\mathcal{L}_{\ell^{(p)}}}

\newcommand*{\cp}[1]{c_{\ell^{(p)}}(#1)}
\newcommand*{\cpfini}[1]{c_{\breve\ell^{(p)}}(#1)}
\newcommand*{\cpfinico}[1]{\hat c_{\breve\ell^{(p)}}(#1)}
\newcommand*{\cpcoe}[1]{c_e(#1)}

\newcommand*{\sqerrlead}{SE_{\ell^{(p)}}}
\newcommand*{\sqerrleadfini}{SE_{\breve\ell^{(p)}}}

\newcommand*{\corr}[2]{\gamma^{#1}\left(j, #2\right)}
\newcommand*{\corrl}[2]{b_S\left(#1,#2,j\right)}
\newcommand*{\corru}[2]{B_S\left(#1,#2,j\right)}

\newcommand*{\rlead}[1][p]{\ell_{\lambda}^{(#1)}}
\newcommand*{\rleadfini}[1][p]{\breve\ell_{\lambda}^{(#1)}}

\newcommand*{\coe}[1]{e_{#1}}
\newcommand*{\coelam}{\coe{\lambda}}
\newcommand*{\coejk}{\coe{j,k}}

\title{
Finite-resolution effects in  $p$-leader multifractal analysis
}

\author{Roberto Leonarduzzi,~\IEEEmembership{Member,~IEEE}, 
Herwig Wendt,~\IEEEmembership{Member,~IEEE}, %
Patrice Abry,~\IEEEmembership{Fellow,~IEEE}, %
Stéphane Jaffard and Clothilde Melot%
\thanks{Manuscript received December 4, 2016; revised March 10, 2017; accepted
  March 17, 2017. Associate editor: Dr. Wei Liu.}
\thanks{R. Leonarduzzi and P. Abry are with  Univ Lyon, Ens de Lyon, Univ Claude Bernard, CNRS,
Laboratoire de Physique, F-69342 Lyon, France ({\tt roberto.leonarduzzi@ens-lyon.fr}, {\tt patrice.abry@ens-lyon.fr}).}%
\thanks{H. Wendt is with IRIT, CNRS UMR 5505, University of Toulouse,  France ({\tt herwig.wendt@irit.fr}).}
\thanks{S. Jaffard is with  Universit\'e Paris Est, Laboratoire d'Analyse et de Math\'ematiques Appliqu\'ees, CNRS UMR 8050, UPEC,  Cr\'eteil, France, ({\tt jaffard@u-pec.fr}).}
\thanks{C. Melot is with the Aix Marseille Univ, CNRS, Centrale Marseille, I2M, Marseille, France ({\tt clothilde.melot@univ-amu.fr}).
.}
\thanks{This work was supported by grant ANR-16-CE33-0020  MultiFracs.
}
\thanks{Digital Object Identifier 10.1109/TSP.2017.2690391
}
}

\begin{document}
\maketitle
\begin{abstract}
\boldmath
Multifractal analysis has become a standard signal processing tool,
for which a promising
new formulation, the $\bm{p}$-leader multifractal formalism, has recently been proposed.
It relies on novel multiscale quantities, the $\bm{p}$-leaders,
defined as local $\ell^p$ norms of sets of wavelet coefficients located at infinitely many fine scales.
Computing such infinite sums from actual finite-resolution data requires  truncations to the finest
available scale, which results in biased $\bm{p}$-leaders and thus in inaccurate estimates of
multifractal properties.
A systematic study of such finite-resolution effects leads to conjecture an explicit and universal
closed-form correction  that permits an accurate estimation of scaling exponents.
This conjecture is formulated from the theoretical study of  a particular class of models for
multifractal processes, the wavelet-based cascades.
The relevance and generality of the proposed conjecture is assessed by numerical simulations
conducted over a large variety of  multifractal processes.
Finally, the relevance of the proposed corrected estimators is demonstrated on the analysis of
heart rate variability data.

\end{abstract}
%


\section{Introduction}
\label{sec:intro}

\paragraph*{Multifractal analysis}
Multifractal analysis has become a standard signal processing tool, widely used and proven relevant in
several different applications, including biomedicine \cite{leonarduzzi2015speech, arneodo2011multi}, finance
\cite{Mandelbrot1999}, geophysics \cite{telesca2011analysis, tessier1993universal}, and
 art investigation \cite{AJWVanGogh2013}, among many others. It amounts to  estimating  the
so-called  \emph{multifractal spectrum} $D(h)$ of a signal or field $X$. 
$D(h)$ quantifies globally and geometrically the local variations of the \emph{regularity} of $X$, measured by the
\emph{regularity exponent} $h$.

\paragraph*{Local regularity}
Traditional formulations  of multifractal analysis rely on the use of the \emph{Hölder exponent} as a
measure of local regularity \cite{Jaffard2004, kantelhardt2002multifractal, muzy1993multifractal}.
However, it has recently been proposed  that multifractal analysis could be based on
\emph{$p$-exponents} instead  \cite{pLeadersPartI2015,pLeadersPartII2015}.
This new formulation presents three key advantages: i) it can be applied to a larger class of
functions or signals $X$ (functions that are locally in $L^p(\RR^d)$ instead of locally bounded);
ii) the variability of the fluctuations of local regularity with $p$ provides additional information on the nature of singularities \cite{Bandt2015,Porque};
 iii) practical estimation methods yield estimates with significantly smaller variance \cite{pLeadersPartII2015}.

\paragraph*{Multifractal formalism}
Estimation of the  multifractal spectrum is conducted in practice following the
so-called
\emph{multifractal formalism} \cite{Jaffard2004,pLeadersPartII2015,Wendt2007d}.
It provides an upper bound for $D(h)$ by analyzing
the scaling behavior \emph{in the limit of fine scales} of appropriate  multiscale quantities, i.e., quantities with a joint
time-scale localization.
Estimates of the corresponding scaling exponents, and thus of $D(h)$, are obtained as linear regressions, over a large range of scales, of time-space averages
of these quantities.

\paragraph*{$p$-leaders}
It has been shown in \cite{pLeadersPartI2015,pLeadersPartII2015} that, when using $p$-exponents, the
multifractal formalism must be based on special multiscale quantities:  the so-called
\emph{$p$-leaders}.
These quantities consist, at a given scale,  of local weighted $\ell^p$ norms of wavelet
coefficients, computed over narrow time neighbourhoods and \emph{over all finer scales}.
The \emph{$p$-leader} multifractal formalism expands and enrich the earlier formulation relying on $\ell^\infty$ norms of wavelet coefficients, the \emph{wavelet leaders} \cite{Jaffard2004,Wendt2007d,wendt_wavelet_2009}.

\paragraph*{Finite-resolution  effect}
When computed from real-world finite-resolution data,  $p$-leaders suffer from two distinct finite-resolution effects. 
The first issue is related to the fact that wavelet coefficients are theoretically defined as
continuous-time inner products, which in practice must be approximated in discrete time.
This subject has already been extensively addressed and shown to have a limited and well-documented impact on the estimation of scaling exponents, cf. e.g.,\cite{abry_initialization_1994,veitch_meaningful_2000}.
The second issue has a much more dramatic impact on estimation quality, and is related to the fact that $p$-leaders, for a given scale, are theoretically defined as sums
\emph{over all the infinitely many finer scales}.
In practice the number of scales is finite, and these sums must necessarily be truncated, yielding a systematic bias
in the actually computed $p$-leaders.
Even worse, this bias is more prominent at  fine scales, which are predominantly
involved in multifractal analysis, thus significantly impairing the estimation of multifractal
parameters.

\paragraph*{Goals, contributions and outline}
The present contribution describes a thorough analysis of finite-resolution effects on
$p$-leader-based multifractal analysis, and proposes a conjecture to practically correct for such effects.
After a short review of the main elements of $p$-leader multifractal analysis in
Section~\ref{sec:MFA}, a conjecture permitting to correct for finite-resolution effects is proposed in Section~\ref{sec:Model}.
It is then
shown theoretically, in Section~\ref{sec:theoretical}, that the proposed conjecture is exact for a special class of multifractal
processes, the wavelet-based cascades.
In Section~\ref{sec:numerical}, the proposed conjecture is further validated by means of numerical simulations
on several multifractal processes of different natures;
$p$-leader estimation performance is also discussed.
Finally, Section~\ref{sec:Illustration} illustrates the relevance of correcting for finite-resolution effects on real data.

\section{Multifractal analysis and $p$-leaders}
\label{sec:MFA}

\subsection{$p$-exponents and multifractal spectrum}
The signal or field to be analyzed is hereafter denoted as $X:\RR^d\to\RR$.
Let $X\in L^p_{loc}(\RR^d)$ for $p\geq 1$. $X$  belongs to $T_{\alpha}^p(x)$, with
$\alpha>-d/p$, if there exist $C,R>0$ and a polynomial $P_{x}$ of degree less than $\alpha$, such that
$\forall a <R$,
$\bigl(\frac{1}{a^d}\int_{B(x,a)} |X(u) - P_x(u)|du \bigr)^{1/p}\leq C a^{\alpha}$, where
$B(x,a)$ is the ball of radius $a$ centered at $x$. The
\emph{$p$-exponent} of $X$ at $x$ is defined as $h_p(x)=\sup\{\alpha\::\:X\in T_{\alpha}^p(x) \}$
\cite{pLeadersPartI2015,CalZyg}. When $p = \infty$, the $p$-exponent
$h_{\infty}(x)$ coincides with the traditional Hölder exponent $h(x)$ \cite{pLeadersPartI2015,JaffMel,Jaffard2004}.
It measures the regularity of  $X$ at $x$: the smaller $h_p(x)$ is, the rougher and
more irregular $X$ is at $x$. Unlike the Hölder exponent, $p$-exponents allow to
measure
\emph{negative regularity}, on condition though that $h_p(x)>-d/p$ \cite{pLeadersPartI2015}.

Multifractal processes are usually defined by the fact that local regularity changes abruptly from one location to another, and a
pointwise estimation of  $h_p(x)$ is therefore of little interest, being itself a highly irregular
function.
Rather, one is interested in a function that quantifies globally the geometrical distribution of
the values $h_p(x)$ takes on:
the multifractal spectrum
$D^{(p)}(h)=\text{dim}_{H}\bigl(\{x\in\RR^d\::\:h_{p}(x)=h\})$, where $\text{dim}_{H}$ denotes the Hausdorff
dimension.
A practical estimate of $D^{(p)}(h)$ requires the use  of multiscale quantities, which we now recall.

\subsection{Wavelet $p$-leaders}

Let $\{\psi^{(i)}(x)\}_{i=1,\ldots,2^d-1}$ denote a family of \emph{mother wavelets}.
These oscillating functions are characterized by a fast decay, good joint time-frequency
localization, and guarantee a \emph{number of vanishing moments} $N_{\psi}\in\NN$, meaning that
$\int x^k\psi^{(i)}(x)dx=0$ for $k=0,1,\cdots,N_{\psi}-1$.
The collection $\{2^{dj/2}\psi^{(i)}(2^jx-k),\;i=1,\cdots,2^d-1,j\in\ZZ,k\in\ZZ^d\}$  of dilated
and translated versions of $\psi^{(i)}$ is an orthonormal basis of
$L^2(\RR^d)$. The discrete wavelet coefficients of $X$ are then defined as:
$\coejk^{(i)}=2^{dj}\int_{\RR^d}X(x)\psi^{(i)}(2^{j}x-k)dx$.
For more details on wavelet transforms, see, e.g., \cite{Mallat1998}.
An $L^1$ normalization
is used in this definition of wavelet coefficients since it is better suited for multifractal analysis.

For simplicity, let $k=(k_1,\dots,k_d)\in\ZZ^d$ and
$\lambda=\lambda_{j,k}=\bigl[2^{-j}k_1,2^{-j}(k_1+1)\bigr)\times\cdots\times\bigl[2^{-j}k_d,2^{-j}(k_d+1)\bigr)$
label \emph{dyadic cubes}.
Each wavelet coefficient
can be associated with
one dyadic cube:
$\coelam^{(i)}=\coejk^{(i)}$.
Let $\lambda(x)$ denote the only cube at scale $j$ that
includes $x$, and   $3\lambda=3\lambda_{j,k}=
\bigl[2^{-j}(k_1-1),2^{-j}(k_1+2)\bigr)\times\cdots\times\bigl[2^{-j}(k_d-1),2^{-j}(k_d+2)\bigr)$
 denote the union of $\lambda$ and its $3^d-1$ neighbours.

Let $p>0$ and $X\in L_{loc}^{p}(\RR^d)$. The \emph{wavelet $p$-leaders} are defined as
\cite{JaffMel,pLeadersPartI2015,pLeadersPartII2015}
\begin{equation}
\label{eq:def_pleaders}
\ell_{j,k}^{(p)} = \ell_{\lambda_{j,k}}^{(p)} := \left(\sum_{\substack{\lambda_{j',k'}'\subset 3\lambda_{j,k}\\ j\leq j'<\infty}}\sum_{i=1}^{2^{d}-1}
\left| \coe{\lambda_{j',k'}}^{(i)}\right|^p 2^{d(j-j')}  \right)^{\frac{1}{p}}
\end{equation}
where $j'\geq j$ is the scale associated with the sub-cube $\lambda'$. Note that the outer sum is
performed over all finer scales $j'\geq j$ and over a narrow spatial neighborhood of $x=2^{-j}k$.

The key property of $p$-leaders is that their decay exactly reproduces the $p$-exponent:
$\ell_{\lambda(x)}^{(p)}\sim 2^{-jh_p(x)}$ when $j\rightarrow\infty$
 \cite{JaffMel,pLeadersPartI2015,pLeadersPartII2015,Jaffard2004}.

 When $p=\infty$, (\ref{eq:def_pleaders}) reduces to the definition of
wavelet leaders, as proposed in \cite{Wendt2007d,wendt_wavelet_2009}.

\subsection{Multifractal formalism}
\label{sec:Formalism}

The \emph{multifractal formalism} permits to estimate $D^{(p)}(h)$ in a
practically feasible and robust way.
It is based on the  so-called \emph{structure
  functions}:
\begin{equation}
\label{eq:def_sf}
\sflead := \frac{1}{n_j}\sum_{k}\left| \ell_{j,k}^{(p)}\right|^q,
\end{equation}
where $n_j$ is the number of coefficients $\ell_{\lambda}^{(p)}$ available at scale $j$.
For multifractal models, $\sflead$ exhibits a  power-law decay, at fine scales,  controlled by the
\emph{scaling exponent} $\zetaq$:
\begin{equation}
\label{eq:def_zeta}
\sflead \sim K_{p,q} 2^{-j\zetaq(q)},\quad j\rightarrow\infty.
\end{equation}
A concave upper-bound for $D^{(p)}$, known as the  \emph{Legendre spectrum} $\legspec$, is provided by the  Legendre transform of
$\zetaq$:
\begin{equation}
\label{eq:def_spectrum}
\legspec(h) :=\inf_{q\in\RR}\bigl(d + qh - \zetaq(q)\bigr) \geq D^{(p)}(h),
\end{equation}
with equality  for numerous multifractal processes, and in particular for the ones used here, cf. \cite{pLeadersPartI2015,pLeadersPartII2015,Bandt2015,Porque}.

\subsection{Log-cumulants}
\label{sec:logcum}
\emph{Log-cumulants} summarize  into a few parameters most of  the relevant information contained in
$\legspec$.
They are defined as  the coefficients of the
Taylor expansion of the scaling function: $\zetaq(q):=\sum_{m\geq
  1}\cp{m}q^m/m!$.
Use of the  Legendre transform also provides an expansion of the
$\legspec$ around its maximum (cf. \cite{wendt_wavelet_2009,pLeadersPartII2015}), further permitting to
interpret the $\cp{m}$:  $\cp{1}$ is the location of the maximum of $\legspec$, $\cp{2}$ is related
to its width, $\cp{3}$ is related to its asymmetry, etc.
By extending calculations in \cite{Castaing1993,Wendt2007d,pLeadersPartII2015}, it is
straightforward to show that the $\cp{m}$ can be computed directly
from the $m$-th order cumulants $\cumjlead{m}$ of $\log \ell_{j,\cdot}^{(p)}$:
\begin{equation}
\label{eq:def_cums}
\cumjlead{m}=C_{\ell^{(p)}, 0} + \cp{m}\log(2^{-j})\quad j\rightarrow\infty.
\end{equation}

\subsection{Practical estimates}
\label{sec:regr}

In practice, $\zetaq(q)$ and $\cp{m}$ are computed by linear regressions, as
$\zetaq(q) = \sum_{j=j_1}^{j_2} \omega_j \log_2\sflead$ and
$\cp{m} = \log_2(e)\sum_{j=j_1}^{j_2} b_j \cumjlead{m}$,
for scales $j$ within the \emph{scaling range} $[j_1,j_2]$, with classical linear regressions weigths
 $b_j$, cf. e.g.,  \cite{wendt_wavelet_2009}.

\subsection{Minimum regularity hypothesis}
\label{sec:min_reg_hyp}
Both $p$-exponents and $p$-leaders are defined only for functions $X\in L^p_{loc}(\RR^{d})$.
It can be easily checked whether data practically satisfy such a property by an \emph{a priori}
analysis of the decay of their \emph{wavelet structure function}
\begin{equation}
\label{eq:def_sf_coef}
\sfcoef := \frac{1}{n_j}\sum_{k}\sum_{i=1}^{2^d-1}\left| \coejk^{(i)}\right|^q, \quad q\geq 0.
\end{equation}
Let $\eta(p)$ denote the \emph{wavelet scaling function}
\begin{equation}
\label{eq:def_eta}
\sfcoef[p] \sim K_{p} 2^{-j\eta(p)}, \quad j\rightarrow\infty.
\end{equation}
It has been shown in \cite{pLeadersPartI2015} that if $\eta(p)>0$, then $X\in L^p_{loc}(\RR^d)$.
It is useful to consider the \emph{critical Lebesgue index} $p_0=\sup(p\::\:\eta(p)>0)$:
$p$-leaders are defined  for $p < p_0$, and when this condition is not met, $p$-leader-based
quantities are not defined theoretically and their practical estimation is thus meaningless
\cite{pLeadersPartI2015,pLeadersPartII2015}.

\section{Finite-resolution effects and estimation}
\label{sec:Model}

\subsection{Finite-resolution scaling behavior}

 Equation (\ref{eq:def_pleaders}) shows that  the computation of
$\ell_{j,k}^{(p)}$ at scale $j$ requires the availability of wavelet coefficients across infinitely many finer
scales $j'$ such that $j\leq j'<\infty$.
However, in practice, only  a
finite-size finite-resolution sampled version of the input data $X$ is available.
Therefore, wavelet coefficients can only be computed for a finite range of scales $\ubar{j}\leq j \leq \bar j$, with $\ubar{j}$ and
$\bar j$ the coarsest and finest scales available.
Thus, the outer sum in (\ref{eq:def_pleaders}) can only be computed for
the \emph{finite} subset of scales $j\leq j'\leq\bar j$, giving rise to \emph{finite-resolution
  $p$-leaders} $\leadfini_{\lambda}$, which suffer from a systematic (under-estimation) bias.

\subsection{Finite-resolution  estimates}
\label{sec:final_estimates}

Let $\sfleadfini$ denote the structure functions  computed from finite-resolution $p$-leaders $\leadfini$.
Motivated by preliminary analyses in \cite{pLeadersPartII2015} and analytical calculations of
wavelet cascades detailed in Section~\ref{sec:theoretical}, we define the following corrected
estimate $\sfleadfinico$:
\begin{align}
\label{eq:sf-nonlinear-corr}
\sfleadfinico &:= \sfleadfini\,\corr{-\frac{q}{p}}{\eta(p)},
\end{align}
\begin{equation}
\label{eq:def_gamma}
\makebox{ with } \, \corr{}{\eta(p)} = \left( \frac{1-2^{-(\bar j- j  +1)\eta(p)}}{1-2^{-\eta(p)}}\right).
\end{equation}
We conjecture that the corrected estimate $\sfleadfinico$ allows to recover the one which would be
obtained from $p$-leaders  computed from infinite-resolution data, i.e.,
\begin{equation}
\label{eq:sf-corr-equivalence}
\sfleadfinico \equiv \sflead.\\
\end{equation}
Equation (\ref{eq:sf-nonlinear-corr}) indicates that scaling in structure functions $S_{\breve\ell^{(p)}}$ computed from
finite-resolution $p$-leaders is corrupted by the nonlinear term
$\gamma$, whose form is conjectured in \eqref{eq:def_gamma}, which can be easily estimated and corrected for.

The following proposition extends Correction (\ref{eq:sf-nonlinear-corr}) to  cumulants.
The proof is sketched in Appendix \ref{app:corr_cums}.

\begin{proposition}
\label{prop:corr_cums}
If and only if  (\ref{eq:sf-nonlinear-corr}) holds, the corrected cumulants $\cumjleadfinico{m}$  relate to
finite- and infinite-resolution cumulants, $\cumjleadfini{m}$ and $\cumjlead{m}$ respectively, as
\begin{align}
\label{eq:c1-nonlinear-corr}
\cumjleadfinico{1} &= \cumjleadfini{1}  - \frac{1}{p}\log\corr{}{\eta(p)},\\
\label{eq:cm-nonlinear-corr}
\cumjleadfinico{m} &= \cumjleadfini{m}\qquad\text{for }m\geq 2,
\end{align}
\begin{equation}
\label{eq:cumj-corr-equivalence}
\cumjlead{m} \equiv \cumjleadfinico{m} \quad\forall m\in\NN^+.
\end{equation}
\end{proposition}

\paragraph*{Remark 1}
The fact that only the scaling of $\cumjleadfini{1}$ is corrupted by finite-resolution effects, while $\cumjleadfini{m}$ for
$m\geq2$ are not, implies that only the
mode $\cp{1}$ of $\legspec$ (i.e., the average regularity) is biased, while the shape (width,
asymmetry, \ldots) is not.
Parameters $\cpfini{m}$ for $m\geq2$, are thus unaffected by finite-resolution effects and benefit from better estimation
performance of $p$-leaders, as detailed in Sec.~\ref{sec:estperf} and also reported in \cite{pLeadersPartII2015}, without
the need of correcting for finite-resolution effects.

\paragraph*{Remark 2}
Because $\corr{}{\eta(p)}$ decays exponentially at coarse scales, $j\rightarrow-\infty$, the
finite-resolution effects become negligible at coarse scales, all the more when $\eta(p)$ is large.

\paragraph*{Remark 3}
When $p\rightarrow\infty$, the proposed correction terms in
(\ref{eq:sf-nonlinear-corr}) and (\ref{eq:c1-nonlinear-corr}) vanish.
Therefore, the conjectured perturbation of scaling at fine scales is not observed for traditional
wavelet leaders (cf. Sections \ref{sec:theoretical} and \ref{sec:numerical} for details).

\vskip2mm\indent%
The following sections show the validity of the proposed corrected estimates, either theoretically
in Sec.~\ref{sec:theoretical} by  analysis of a special class of multifractal processes, the wavelet
cascades, or empirically in Sec.~\ref{sec:numerical} by means of numerical simulations conducted
over several   multifractal processes  different in nature.

\section{Theoretical results}
\label{sec:theoretical}

In this section, finite-resolution effects are investigated theoretically on functions defined directly by
wavelet coefficients, i.e.,  1D and 2D deterministic or random wavelet cascades,
for which $S_{\coe{}}$ and $C_{\coe{}}$ can be regarded as the exact scaling quantities.

For ease of exposition, this section makes use of the  \emph{restricted $p$-leaders} $\rlead$,
defined by replacing $3\lambda$ with $\lambda$ in (\ref{eq:def_pleaders}):
\(
\ell_{j,k}^{(p)} = \ell_\lambda^{(p)} = \left(\sum_{\lambda'\subset \lambda}\sum_{i=1}^{2^{d}-1}
|c_{\lambda'}^{(i)}|^p 2^{d(j-j')}  \right)^{1/p}
\).
It has been shown that structure functions computed with restricted $p$-leaders and $p$-leaders as in \eqref{eq:def_pleaders} are
equivalent (cf.~\cite{Bergou}).

\subsection{Deterministic Binomial Wavelet Cascade}
\label{sec:dbwc}

\subsubsection{Construction}
Inspired by \cite{decoster2000}, we propose a model for 2D Deterministic Binomial Wavelet Cascade (DBWC), whose wavelet
coefficients are defined as follows:
\begin{equation}
\label{eq:def_dbwc2d}
\left\{
\begin{array}{ll}
d_{0,1,1} &= 1\\
d_{j,2k_1,2k_2} &= w_0\, d_{j-1,k_1,k_2}\\
d_{j,2k_1+1,2k_2} &= w_1\, d_{j-1,k_1,k_2}\\
d_{j,2k_1,2k_2+1} &= w_2\, d_{j-1,k_1,k_2}\\
d_{j,2k_1+1,2k_2+1} &= w_3\, d_{j-1,k_1,k_2}\\
\coe{j,k_1,k_2}^{(i)}&=\alpha_i\, d_{j,k_1,k_2}
\end{array}
\right.
\end{equation}
with weights $w_i$ being deterministic constants controlling  multifractal properties, and $\alpha=(\alpha_i)$ controlling anisotropy.
It can be shown that the wavelet scaling function $\eta$ reads (cf. Appendix~\ref{app:dbwc_eta}),
\begin{equation}
\label{eq:eta_dbwc_2d}
\eta(q) = 2 -\log_2\sum_{m=0}^3w_m^q,\quad\text{for }q>0
\end{equation}
and that  DBWC satisfies $\zetaq(q)=\eta(q)$, $q>0$,
while anisotropy 
has no impact on the scaling properties.

\subsubsection{$p$-leader analysis}

Finite-resolution effects for 2D DBWC are described by the following proposition.
\begin{proposition}
\label{prop:dbwc_est}
The finite-resolution $p$-leader structure functions  of a 2D DBWC as in (\ref{eq:def_dbwc2d}) are given by
\begin{equation}
\label{eq:sf-dwbc}
\sfleadfini = \Vert\alpha\Vert_p \;\sfcoef\; \corr{\frac{q}{p}}{\eta(p)},
\end{equation}
where the function $\gamma$ is defined in \eqref{eq:def_gamma}.
\end{proposition}

The proof, cf. Appendix~\ref{app:dbwc_prop}, relies on the multiplicative structure of wavelet
coefficients.
Comparing \eqref{eq:sf-dwbc}  with (\ref{eq:sf-nonlinear-corr}) shows the relevance of the proposed correction.
Similar computations for 1D DBWC lead to identical conclusions,
with notably the same correction function $\gamma$ (cf. \cite{pLeadersPartII2015}).

\subsection{Multiplicative Random Wavelet Series (MRWS)}
\label{sec:rws}

\subsubsection{Construction}
Random wavelet series (RWS), originally introduced in \cite{AJ02}, are a general framework for constructing multifractal functions from their wavelet expansion.
They are built by assigning to each wavelet coefficient an independent realization of a random
variable. Here, we will consider the specific case of multiplicative RWS. Let $\{W\}^{(j)}$ denote the product of $j$ independent copies of the continuous
positive random variable $W$. Let $e_{0,1}=1$.
Then, the $2^{j}$ coefficients at scale $j>0$ are built as $\coejk\stackrel{\mathcal{L}}{=}\{W\}^{(j)}$.
The wavelet scaling function $\eta$, for $q>0$, reads $\eta(q)=-\log_2\EE{W^q}$
under suitable assumptions on the tail of $W$ \cite{Arneodo1998rwc},
and MRWS satisfy $\zetaq(q)=\eta(q)$, for $q>0$.

\subsubsection{$p$-leader analysis}

First, we analyze the behavior of $\sfleadfini$ for $q$ a multiple of $p$.

\begin{proposition}
\label{prop:rws-sf-nn}
Let $q=np$, $n\in\NN$.
The MRWS finite-resolution $p$-leader structure function $\sfleadfini $ satisfies:
\begin{multline}
\label{eq:rws-sf-nn}
\corrl{n}{p} \leq
\frac{\sfleadfini }{\sfcoef \,\corr{n}{\eta(p)}}
\leq \corru{n}{p},
\end{multline}
with
\begin{align}
\label{eq:def-corr-sf-low}
\corrl{n}{p}&=2^{-j(n\eta(p)-\eta(np))}, \\
\label{eq:def-corr-sf-up}
\corru{n}{p}&=\frac{\corr{n}{\eta(np)/n}}{\corr{n}{\eta(p)}},
\end{align}
and where the function $\gamma$ is defined in \eqref{eq:def_gamma}.
\end{proposition}

The proof, cf. Appendix~\ref{app:rws-sf-nn}, relies on the multiplicative structure of wavelet
coefficients and the concavity of the scaling function.
Proposition \ref{prop:rws-sf-nn} can be extended to all positive values of $q$, as
in the following consequence.

\begin{consequence}
\label{prop:rws-sf-rr}
Assuming that $\sfleadfini[np]$ and $\sfcoef[np]$ are smooth enough as functions of $n$,
then Prop.~\ref{prop:rws-sf-nn} also holds for $n\in\RR^+$.
\end{consequence}

The proof, cf. Appendix~\ref{app:rws-sf-rr}, is based on an interpolation argument.

\subsubsection{Remarks}
\begin{figure}
\begin{tabular}{*{2}{m{40mm}}}
\begin{overpic}[width=\linewidth]{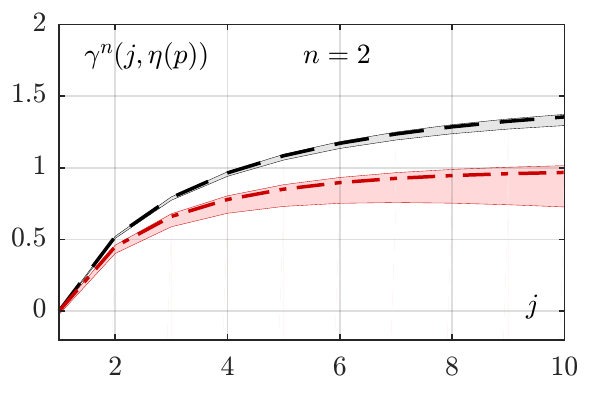}
  \put(23,12){\includegraphics[scale=0.4]{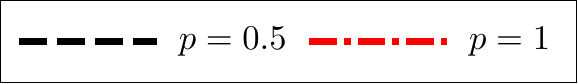}}
\end{overpic}
&
\includegraphics[width=\linewidth]{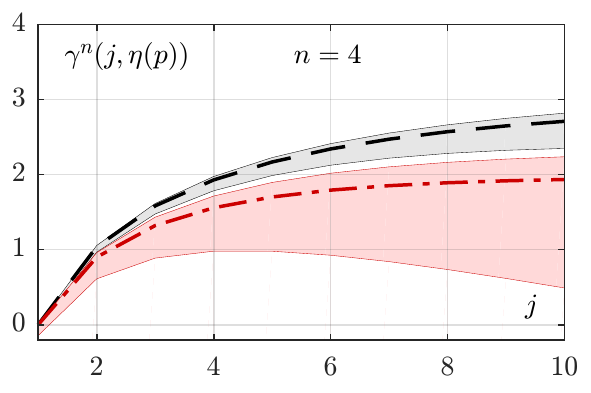}
\end{tabular}
\caption{\textbf{MRWS structure functions.} Proposed correction $\corr{n}{\eta(p)}$, for two values
  of $n$ (left and right panels), and $p=0.5$ (black dashed line) and $p=1$ (red dash-dotted line).
  The solid lines that delimit the shaded areas represent the bounds $\corrl{n}{p}$  and
  $\corru{n}{p}$, which converge to each other for small $j$ and $p$.}
\label{fig:rws-bounds}
\end{figure}

Proposition \ref{prop:rws-sf-nn} shows that, for MRWS, we
are only able to produce bounds for the deviation from exact scaling induced by finite-resolution
effects.
However, the bounds $b_S$ and $B_S$ tend to coincide for small values of $p$, $q=np$ and $j$, as
illustrated in Fig.~\ref{fig:rws-bounds}.
Thus, the proposed corrections can be assumed to be asymptotically exact in those situations when the
finite-resolution effects are the strongest (small $p$ and fine scales $j\rightarrow\bar{j}$).
Also, the lower and upper bounds coincide when $\eta$ is a linear function, indicating that the
proposed corrections are exact for monofractal MRWS.

\subsection{Random Wavelet Cascades (RWC)}
\label{sec:rwc}

\subsubsection{Construction}
The MRWS analyzed in the previous section have independent wavelet coefficients.
We now consider a related process with strongly correlated wavelet coefficients:
the random wavelet cascades (RWC) \cite{Arneodo1998rwc}.
Let $\coejk=1$, and
let $W_l$, $W_r$ and $W$ denote iid  positive random variables.
Wavelet coefficients at scales  $j=1,2,\cdots,\bar j$ are built from
coefficients at scale $j-1$ by the iterative procedure $\coe{j,2k}=W_l\coe{j-1,k}$,
$\coe{j,2k+1}=W_r\coe{j-1,k}$.
The wavelet scaling function $\eta$, defined for $q>0$, is also shown to read
$\eta(q)=-\log_2\EE{W^q}$ under suitable assumptions on the tail of $W$ \cite{Arneodo1998rwc},
and RWC satisfies $\zetaq(q)=\eta(q)$ for $q>0$.

\subsubsection{$p$-leader analysis}

The complicated correlation structure precludes the  computation of the structure functions for
an arbitrary $q$.
Thus, we restrict calculations to $q=p$ and $q=2p$.

\begin{proposition}
\label{prop:rwc_sf_p}
For $q=p$, the $p$-leader structure function of a RWC reads
\begin{equation}
\label{eq:rwc-sf-p}
\sfleadfini[p] = \sfcoef[p]  \,\corr{}{\eta(p)}
\end{equation}
where the function $\gamma$ is defined in \eqref{eq:def_gamma}.
\end{proposition}

\begin{proposition}
\label{prop:rwc_sf_2p}
For $q=2p$, the $p$-leader structure function of a RWC reads
\begin{equation}
\label{eq:rwc-sf-2p}
\sfleadfini[2p] = \sfcoef[2p]   \corr{2}{\eta(p)}\, f(j,p),
\end{equation}
where
\begin{equation}
\label{eq:rwc-extra-term}
f(j,p) =  \frac{1}{2-2^{-\mu(p)+1}} \left[1   - 2^{\mu(p)}\frac{\corr{}{\mu(p)}}{\corr{}{\eta(p)}}\right],
\end{equation}
with $\mu(p)=\eta(2p)-\eta(p)+1$, and $\gamma$ is defined in \eqref{eq:def_gamma}.
\end{proposition}

\begin{figure}
\centering
\begin{tabular}{cc}
\raisebox{-0.5\height}{\includegraphics[width=0.6\linewidth]{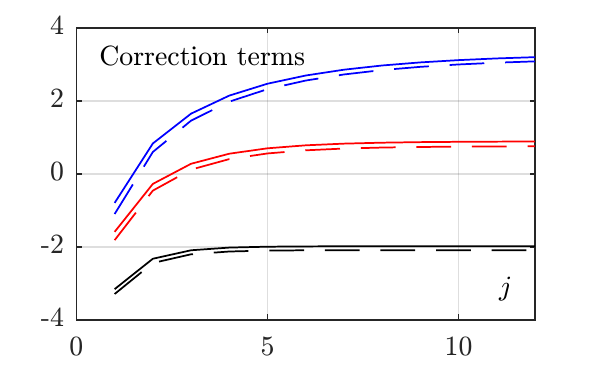}}&
\raisebox{-0.5\height}{\hspace*{-7mm}\includegraphics[scale=0.7]{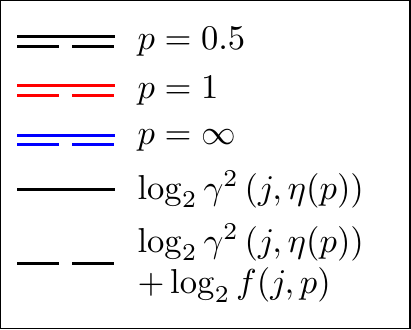}}\\
\end{tabular}
\caption{\textbf{Random wavelet cascades, \bm{$q=2p$}.} Correction terms
  $\log_2\corr{2}{\eta(p)}$ (solid)  and $\log_2\corr{2}{\eta(p)} f(j,p)$
    (dashed), for several values of $p$ (colors). The difference is negligible.}
\label{fig:rwc_corr_comp}
\end{figure}

Proofs are given in Appendix \ref{app:rwc-sf-2p}.
Proposition \ref{prop:rwc_sf_2p} shows that, in the presence of correlations, the finite-resolution effect
theoretically differs from the correction conjectured in (\ref{eq:sf-nonlinear-corr}) by the higher-order term
$f(j,p)$.
However, extensive numerical simulations indicate that this term has negligible effect.
This is illustrated in Fig.~\ref{fig:rwc_corr_comp}, using $\eta(p)=c(1)p+c(2)^2/2$, a typical
example for many processes (here with
$c(1)=0.8$ and $c(2)=-0.08$).

\section{Empirical assessment} 
\label{sec:numerical}

\subsection{Multifractal processes}
In this section, we investigate the level of validity of the proposed corrected estimates of
Sec. \ref{sec:final_estimates} in a general setting, using a  representative panel of  multifractal processes.

\subsubsection{Fractional brownian motion (fBm)}
fBm is defined as the integral of a Gaussian noise with a kernel that
defines its covariance structure \cite{mandelbrot_fractional_1968,Samorodnitsky1994}, fully controlled by
the \emph{Hurst parameter} $H$.
fBm is monofractal, which means that its  $D^{(p)}$ collapses to a single point.

\subsubsection{Multifractal random walk (MRW)}
MRW is defined from two independent Gaussian processes, with a specific covariance structure
chosen to mimic that of multiplicative cascades  \cite{bacry2001}.
A 2D extension has been  proposed in \cite{chevillard_stochastic_2010}.
It has a  parabolic $D^{(p)}$ controlled by two parameters, $H$ and $\lambda$.
Expressions for its multifractal spectrum and $p_0$ are provided in
\cite{bacry2001,chevillard_stochastic_2010,pLeadersPartII2015}.

\subsubsection{$\alpha$-stable Lévy process}
An $\alpha$-stable Lévy process is defined  as a selfsimilar process with independent stationary increments
 \cite{Samorodnitsky1994}.
It has a linear $D^{(p)}$, controlled by the selfsimilarity exponent $\alpha$.
Expressions for the multifractal spectrum and $p_0$ are given in \cite{JaffardLevy,LeonarduzziGretsi2015}.

\paragraph*{Critical Lebesgue index}
The  critical Lebesgue index $p_0$ is always $\infty$ for the considered processes.
For MRW and $\alpha$-stable Lévy process, we will also analyze their fractional derivatives,
which have a finite $p_0$ tuned by the differentiation order (cf.~\cite{pLeadersPartII2015,LeonarduzziGretsi2015}).

\subsection{Simulation setup}
$N_{MC}=100$ realizations of each multifractal process are analyzed, of size $N=2^{19}$ for 1D
processes and $N_1\times N_2=2^{10}\times2^{10}$ for 2D processes.
In all cases, averages
over the $N_{MC}$ realizations are reported.
The synthesis parameters were set to
$H=0.7$ for fBm,
 $H=0.84$ and $\lambda=\sqrt{0.08}$ for MRW (both 1D and 2D), and
$\alpha=0.8$ for Lévy process.

Wavelet analysis is performed using a Daubechies wavelet with $N_{\psi}=3$ vanishing moments.
$p$-leaders are computed
for $p \in \{1/4,1/2, 1, 2, 5, \infty\}$, and the
 convention that the finest available scale is $\bar j=1$.
Scaling exponents and log-cumulants  are computed using weighted linear
regressions \cite{Wendt2007d}.

\subsection{Logscale diagrams}

\subsubsection{Impact of the proposed correction for $C(1,j)$}
\begin{figure}[t]
\centering

\providecommand{\myinc}[1]{
  \includegraphics[width=1\linewidth,trim=4 5 2 2,clip=true]{#1}
}
\begin{tabular}{*{2}{m{4cm}}}
\myinc{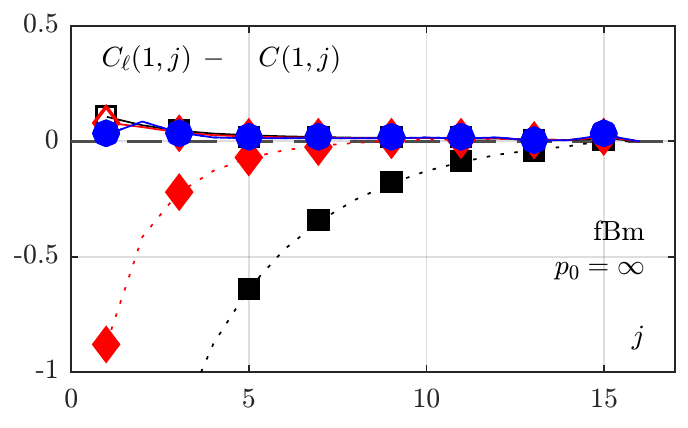} &
\hspace*{5mm}\includegraphics[scale=0.6]{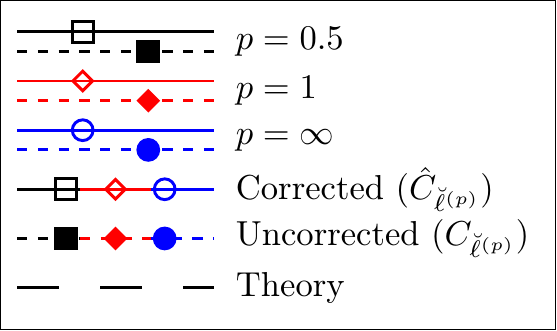} \\
\myinc{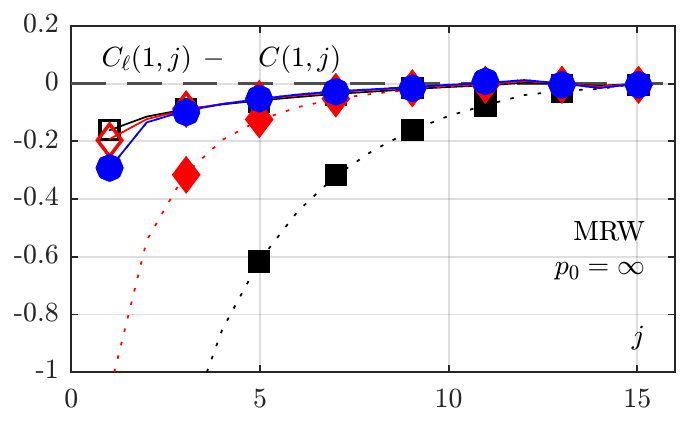} &
\myinc{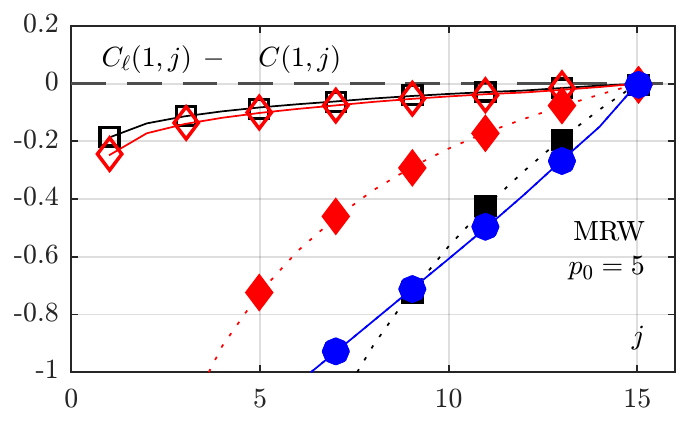}\\
\myinc{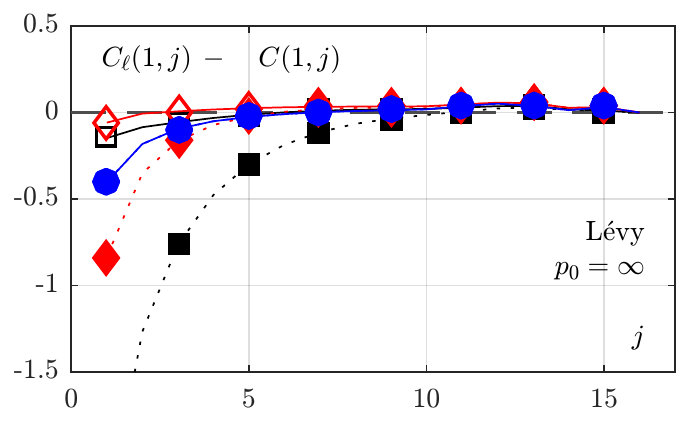} &
\myinc{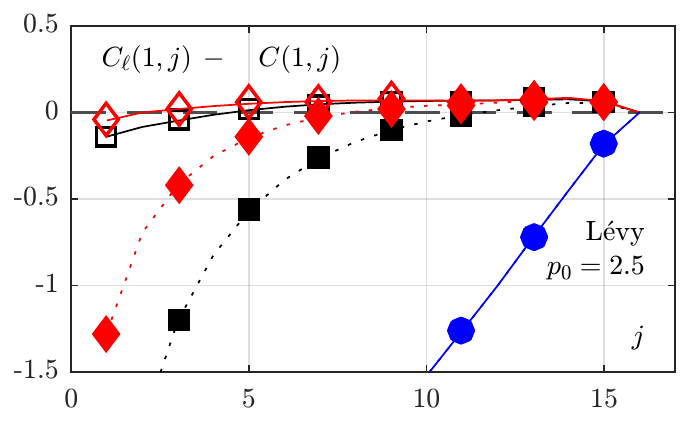} \\
\myinc{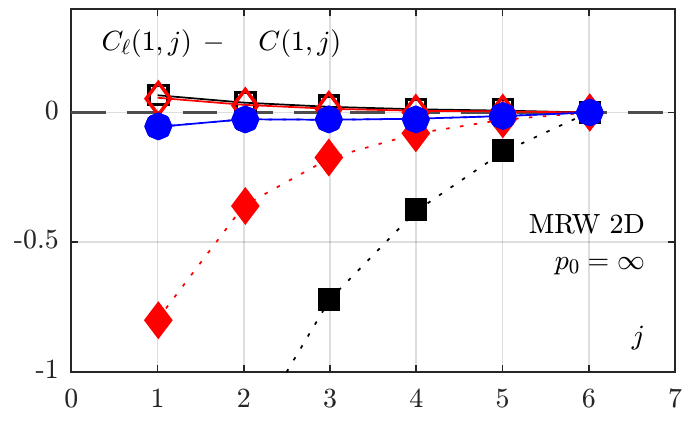} &
\myinc{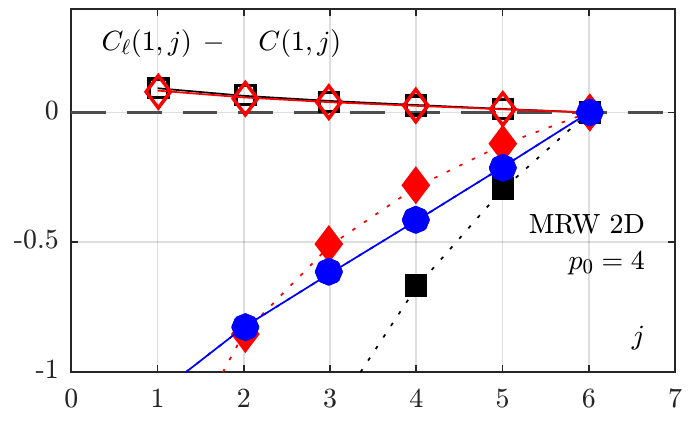}
\end{tabular}
\caption{\textbf{Logscale diagrams: impact of correction.}
  Logscale diagrams  for different processes (rows) and values of $p_0$ (columns).
  Solid lines with empty markers represent  corrected $\cumjleadfinico{1}$, while dashed lines with
  solid markers represent uncorrected $\cumjleadfini{1}$.
  Marker-styles and colors indicate different values of $p$.
}
\label{fig:loglogs}
\end{figure}

We begin by analyzing qualitatively and quantitatively how the proposed correction enables  to restore the correct scaling behavior for $C(1,j)$.
Fig.~\ref{fig:loglogs} superimposes cumulants with correction $\cumjleadfinico{1} - C(1,j)$ (solid
lines, empty markers) and without correction $\cumjleadfini{1} - C(1,j)$  (dotted lines, full markers), for several values of  $p$ and $p_0$.
The subtraction of the true scaling $C(1,j)$ is intended to ease comparisons since departures from perfect estimation thus materialize as departures from $0$.
Fig.~\ref{fig:loglogs} strikingly shows that uncorrected $\cumjleadfini{1}$ present significant
departures from the theoretical scaling, and clearly depart one from another for different values of
$p$.
To the contrary,  corrected $\cumjleadfinico{1}$ show very mild departures from the theoretical
scaling, and additionally they all coincide.
These are very satisfactory outcomes as it is known theoretically that for all processes analyzed here the multifractal spectra  $D^{(p)}(h)$ (and hence $C(1,j)$) do not depend on $p$.
These observations suggest that the conjectured correction (\ref{eq:c1-nonlinear-corr}) is
valid and effective for a much larger class of processes  than those studied in Section
\ref{sec:theoretical}.
%
Uncorrected $\cumjleadfini{1}$ yield departures from theoretical behavior that are larger for small $p$ as well
as for $p_{0}<\infty$,
which is consistent with the fact that $\eta(p)$ is smaller for small values of $p$ and of $p_0$.
\paragraph*{Remark}
Despite the fact that wavelet leaders ($p=\infty$) are not defined for
$p_0<\infty$ \cite{pLeadersPartI2015}, they can still be computed in practice. However---as shown in
Fig.~\ref{fig:loglogs}(right)--- these practical estimates are affected by a strong bias, which is
explicitly accounted for in \cite{LeonarduzziGretsi2015,pLeadersPartII2015}.

\subsubsection{Quantitative assessment}
To further assess the relevance of the correction on the logscale diagrams we propose to quantify
the deviations from true scaling by the  squared error
$\sqerrlead=\langle \sum_{j=\bar j}^{\ubar{j}} ( \cumjleadfinico{1} - C(1,j) )^2\rangle_N$,
where $\langle\cdot\rangle_N$ stands for the average over $N$ independent realizations, and compare
them to the deviations $\sqerrleadfini$ yielded when no correction is used, i.e., by
$\cumjleadfini{1}$.
Because only the scaling behavior is of interest here, the influence of the intersect $C^{(p,0)}(1)$
is removed, by simple substraction so that
$\cumjlead[\underline{j}]{1}-C(1,\underline{j})=0$ (and similarly for $\cumjleadfini{1}$).
Table \ref{table:loglog_error} reports results in terms of $\log_{10} (\sqerrlead / \sqerrleadfini)$ for several $p_0$ and $p$.
The fact that most entries in Table \ref{table:loglog_error} are positive confirms that the use of
the correction lessens the difference with
the correct scaling behavior. 
Notably, for small $p$ and $p_0$, the corrected estimator $\cumjleadfinico{1}$ improves the
squared error by 3 orders of magnitude over the uncorrected $\cumjleadfini{1}$.
For larger $p$ and
$p_0$, where the impact of correction appears to be less significant, the scaling of
$\cumjleadfini{1}$ is actually already close to the theoretical one; correction is thus less
needed.

\begin{table}[t]
\centering
\caption{\textbf{Relative squared error in departures from scaling.} $\log_{10} (\sqerrlead / \sqerrleadfini)$,
where $\sqerrlead$ and $\sqerrleadfini$ are the squared errors quantified the departures of $\cumjleadfinico{1}$ and $\cumjleadfini{1}$ from the exact theoretical scaling.}
\label{table:loglog_error}
\begin{tabular}{ *{6}{c} }
& & \multicolumn{4}{c}{$p$} \\ \cline{3-6}
& $p_0$ &  $0.5$ & $1$ & $2$ & $5$ \\
\hline
fBm
& $\infty$ & $ 2.81$ & $ 1.72$ & $0.28$ & $-0.03$ \\
\hline
\multirow{3}{*}{\rotatebox[origin=c]{90}{MRW}}
& $ 1.3$   &$ 2.39$ &$ 1.76$ &$ 1.25$  &$ 1.51$ \\
& $  2.5$   &$ 2.97$ &$ 2.15$ &$ 1.42$ &$ 1.56$ \\
& $    5$   &$ 2.71$ &$ 1.82$ &$ 0.99$ &$ 0.61$ \\
& $  \infty$&$ 2.35$ &$ 1.29$ &$0.46$  &$0.05$ \\
\hline
\multirow{3}{*}{\rotatebox[origin=c]{90}{Lévy}}
& $ 1.3$   &$ 3.66$ &$ 2.53$ &$ 2.40$  &$ 2.77$ \\
& $  2.5$   &$ 3.09$ &$ 1.85$ &$ 0.80$ &$ 1.54$ \\
& $    5$   &$ 3.01$ &$ 1.89$ &$ 1.00$ &$-0.13$ \\
& $  \infty$&$ 2.71$ &$ 1.62$ &$ 1.30$ &$0.93$ \\
\hline
\multirow{3}{*}{\rotatebox[origin=c]{90}{MRW 2D}}
& $  1.5$   &$ 3.31$ &$ 2.76$ &$ 2.39$  &$ 2.70$ \\
& $ 2.8$   &$ 2.74$ &$ 2.13$ &$  1.6$  &$ 1.41$ \\
& $ 5.3$   &$ 3.05$ &$ 2.42$ &$ 1.88$  &$ 2.02$ \\
& $  \infty$&$ 2.93$ &$ 2.13$ &$ 1.25$ &$0.55$ \\
\hline
\end{tabular}
\end{table}

\subsubsection{Logscale diagrams for $C(m,j)$, $m\geq 2$}
Fig.~\ref{fig:loglogsc2} provides examples of $\cumjleadfini{2}$ for MRW, for two different
critical Lebesgue indices $p_0$, and different $p$s.
Fig.~\ref{fig:loglogsc2} clearly shows that $\cumjleadfinico{2}$, for all $p$, reproduce the expected theoretical scaling, as functions of scales $j$, independently of $p_0$, and that, as expected, the $\cumjleadfinico{2}$ superimpose for all $p$.
This confirms numerically that  no finite-resolution effects are observed on the higher-order
cumulants $\cumjleadfini{m}$, $m\geq2$, and thus
no correction is needed, cf. \eqref{eq:cm-nonlinear-corr}.

\subsubsection{Structure functions $\sfleadfini$ and $\sflead$}
Since the scaling of structure functions can be directly translated into the scaling of cumulants
(cf., Section \ref{sec:logcum} and Proposition \ref{prop:corr_cums}), the relevance of the correction  for the structure
function $\sflead$ is directly determined by the relevance and accuracy of the correction for the
first cumulant, $\cumjleadfinico{1}$, which has been extensively assessed above. Therefore, the above
results and conclusions for $\cumjleadfinico{1}$ directly apply to structure functions, and are  not reproduced
or further discussed here.

\vskip2mm

Overall, these results unambiguously indicate that the conjectured
corrections (\ref{eq:sf-nonlinear-corr}-\ref{eq:def_gamma}) generically and robustly enable to
remove the finite-resolution bias from cumulants  and structure functions, and to restore their expected
scaling behavior.

\providecommand{\myinc}[1]{
\includegraphics[width=1\linewidth,trim=3 4 2 2,clip=true]{#1}
}
\begin{figure}[t]
\centering
\begin{tabular}{*{2}{m{4cm}}}
\multicolumn{2}{c}{\includegraphics[scale=0.5]{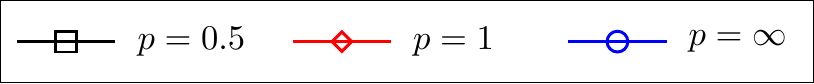}}\\
\myinc{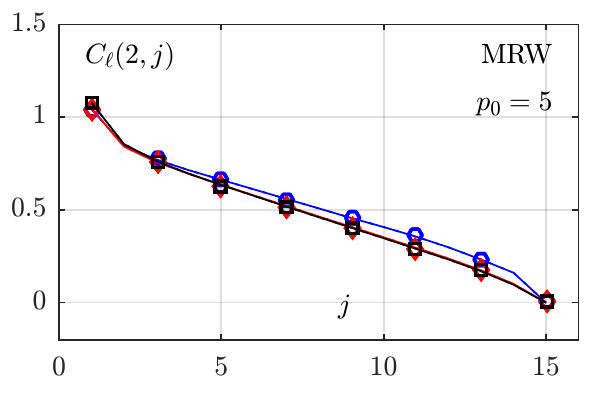}&
\myinc{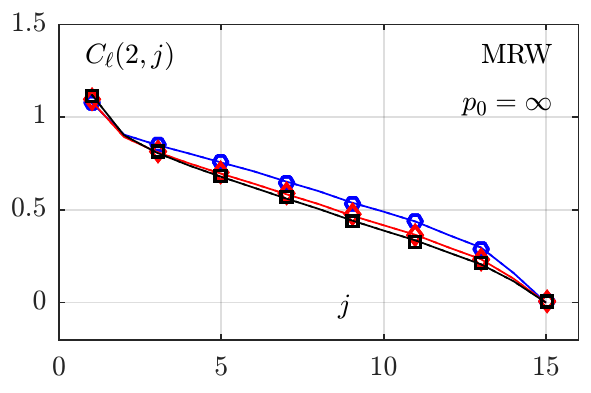}
\end{tabular}
\caption{\textbf{\boldmath Uncorrected logscale diagrams $\cumjleadfini{2}$}, for MRW, for different values of $p$ and $p_0$.
}
\label{fig:loglogsc2}
\end{figure}

\subsection{Estimation performance for scaling parameters}
\label{sec:estperf}

\subsubsection{Estimation of $c(1)$}
Estimation of scaling parameters requires the selection of a range of scales where the linear
regression is performed.
As suggested in Fig.~\ref{fig:loglogs}, the impact of  finite-resolution effects
on the bias of scaling-parameter estimates could be reduced by performing linear regressions
at sufficiently coarse scales, yet at the price of a significant increase of the corresponding estimation variance.
To quantify this, we set the upper limit of the scaling range $j_2$ to the coarsest available scale, and evaluate estimation
performance for linear regressions conducted from all possible lower limits $j_1$, with both
corrected and uncorrected cumulants  $\cumjleadfinico{1}$ and  $\cumjleadfini{1}$.

\begin{figure}
\centering
\noindent\begin{tabular}{ll}
\raisebox{-0.5\height}{\includegraphics[width=0.6\linewidth]{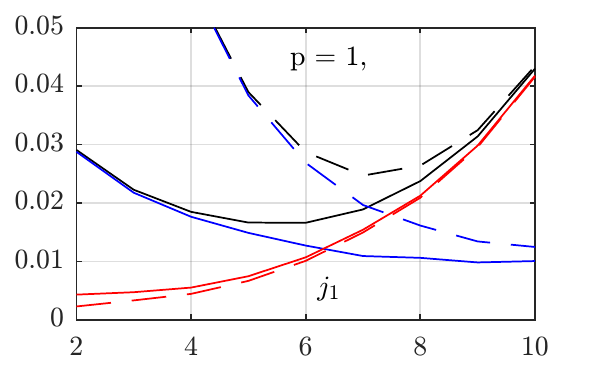}} &
\raisebox{-0.5\height}{\hspace*{-8mm}\includegraphics[scale=0.7]{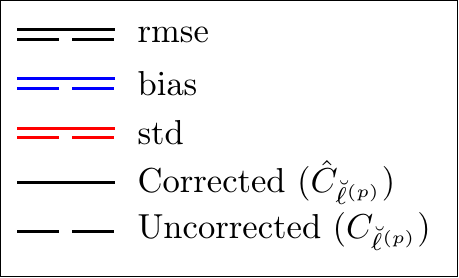}}
\end{tabular}
\caption{\textbf{Estimation performance (bias, std, rmse) for estimates $\overline{c}_1^{(p)}$ and
    $\hat{c}_1^{(p)}$ as functions of lower scale $j_1$.} Solid lines represent corrected
  $\cumjleadfinico{1}$, while dashed lines represent uncorrected $\cumjleadfinico{1}$. The upper
  scale $j_2$ was set to the largest available scale.}
\label{fig:perf_j1}
\end{figure}

Estimation performances for $c(1)$ as functions of $j_1$ are compared in
Fig.~\ref{fig:perf_j1}, for MRW, in terms of bias, standard deviation (std) and root mean squared
error (rmse).
Benefits of the proposed correction on estimation performance are striking. First,
bias is significantly reduced for $\cpfinico{1}$ as compared to that  of $\cpfini{1}$, which
is subject to a dramatic blow-up for small values of $j_1$. Second,
correction for the bias does not alter the std.
Consequently, the smallest rmse for $\cpfinico{1}$
 is achieved at  $j_1=5$  while only at scale $j_1=7$ for
$\cpfini{1}$, i.e., the correction enables the \emph{use of finer scales} in linear
regression;
moreover, the optimal $rmse$ is smaller for  $\cpfinico{1}$ than for $\cpfini{1}$.

%
%
%

To further quantify the decrease in rmse   and in usable fine scales yielded by the conjectured correction, 
Table \ref{tab:optimal_mse} (left panel) reports the relative optimal
rmse (RORMSE) for $\cpfinico{1}$ and $\cpfini{1}$, defined as
$\textnormal{RORMSE}= \min_{j_1}\textnormal{rmse}(j_1;\cpfini{1}) /
\min_{j_1}\textnormal{rmse}(j_1;\cpfinico{1})$.
Table \ref{tab:optimal_mse} clearly demonstrates that using the correction (\ref{eq:c1-nonlinear-corr})
can yield considerable reductions of rmse values, by up to one order of magnitude.
The gains in rmse are smaller for large $p$ and $p_0$, as can be expected from the fact that in these
cases $\eta(p)$ takes on large values and finite-resolution effects are hence negligible.
In Table \ref{tab:optimal_mse} (right panel), the choices of fine scale $j_1$ that lead to best rmse
values, denoted optimal lower cutoff (OLC) and defined as
$\textnormal{OLC}(c(1))=\arg\min_{j_1}\textnormal{rmse}(j_1;c(1))$, are compared for $\cpfini{1}$
and $\cpfinico{1}$.
The OLC values indicate that the conjectured correction indeed permits the use of several additional
finer scales for scaling parameter estimation, thus explaining the origin in the reduction of RORMSE.

\subsubsection{Estimation of $c_m$, $m\geq 2$}

Table \ref{tab:cm} reports the rmse for estimates  $\cpfinico{m}$, for
$m=1,2,3$, for MRW and different $p$ and $p_0$, and shows that:
i) the rmse for $\cpfinico{m}$  always decreases with $p$;
ii) while the rmse
of $\cpfinico{1}$ and $\cpcoe{1}$ are similar, there is a large difference  for $m\geq 2$,
illustrating the inability of wavelet coefficients to estimate higher-order log-cumulants,
and iii) while decreasing $p_0$ increases the rmse of $\cpfinico{1}$ when $p>p_0$, it is not the
case for higher-order log-cumulants (cf. \cite{LeonarduzziGretsi2015} for details).
These results clearly indicate that $p$-leader multifractal analysis with small $p$ always yields
the best estimation performance (see also \cite{pLeadersPartII2015}).
\vskip2mm\indent

\subsubsection{Importance of accounting for finite-resolution effects}

This section demonstrated, first, that corrections
(\ref{eq:sf-nonlinear-corr}) and (\ref{eq:c1-nonlinear-corr}) are robust and valid for large classes of processes
and, second, that they permit a dramatic improvement in the accuracy of
$p$-leader-based scaling analysis by: i) significantly reducing estimation bias, whatever
$j_1$; ii) allowing to use several additional fine scales; iii) reducing rmse for the estimation of $c(1)$ by up to an order of magnitude.

\begin{table}[t]
\centering
\caption{\textbf{Optimal rmse and lower cutoff.}
Left panel: relative optimal rmse RORMSE; larger values indicate a larger gain due to the correction term.
Right panel: optimal lower cutoffs $\textnormal{OLC}(c(1))=\arg\min_{j_1}\textnormal{rmse}(j_1;c_{1})$  for $\cpfini{1}$ and $\cpfinico{1}$.
 }
\label{tab:optimal_mse}
\begin{tabular}{*{2}{c}||*{3}{c}||*{3}{c}}
& &
 \multicolumn{3}{c||}{\boldmath$\log_{10}($\textbf{RORMSE}$)$}
  &
 \multicolumn{3}{c}{\boldmath{$\frac{  \text{\bf OLC}(\cpfini{1})   }{   \text{\bf OLC}(\cpfinico{1})}$}}
\\ \hline
& &
 \multicolumn{3}{c||}{\boldmath{$p$}} & \multicolumn{3}{c}{\boldmath{$p$}}\\
\cline{3-8}
& \boldmath{$p_0$} & \boldmath{$0.5$} & \boldmath{$2$} & \boldmath{$5$} & \boldmath{$0.5$} & \boldmath{$2$} & \boldmath{$5$} \\
\hline
fBm
 & \boldmath{$ \infty$} & $0.98$   &   $-0.11$   &   $0.00$ & $9 \,/\, 5$    &   $2 \,/\, 3$   &   $3 \,/\, 3$ \\
\hline
\multirow{4}{*}{\rotatebox[origin=c]{90}{MRW}}
 & \boldmath{$ 1.3$}     &   $0.85$   &   $0.60$   &   $0.62$ & $9 \,/\, 5$   &   $9 \,/\, 5$   &   $8 \,/\, 5$ \\
 & \boldmath{$ 2.5$}     &   $0.81$   &   $0.52$   &   $0.50$ & $9 \,/\, 6$   &   $9 \,/\, 6$   &   $8 \,/\, 5$ \\
 & \boldmath{$   5$}     &   $0.63$   &   $0.30$   &   $0.17$ & $9 \,/\, 5$   &   $8 \,/\, 6$   &   $8 \,/\, 7$ \\
 & \boldmath{$ \infty$}  &   $0.17$   &   $0.02$   &   $0.00$ & $9 \,/\, 4$   &   $5 \,/\, 4$   &   $4 \,/\, 4$ \\
\hline
\multirow{4}{*}{\rotatebox[origin=c]{90}{Lévy}}
 & \boldmath{$ 1.2$}    &   $0.59$   &   $0.81$   &   $0.83$  & $9 \,/\, 3$   &   $7 \,/\, 1$   &   $5 \,/\, 1$ \\
 & \boldmath{$ 2.5$}    &   $0.30$   &   $0.43$   &   $0.46$  & $9 \,/\, 3$   &   $7 \,/\, 1$   &   $3 \,/\, 1$ \\
 & \boldmath{$   5$}    &   $0.19$   &   $0.13$   &   $-0.17$ & $9 \,/\, 5$   &   $5 \,/\, 2$   &   $5 \,/\, 1$ \\
 & \boldmath{$ \infty$} &   $0.18$   &   $0.00$   &   $-0.51$ & $9 \,/\, 4$   &   $3 \,/\, 2$   &   $3 \,/\, 1$ \\
\hline
\multirow{4}{*}{\rotatebox[origin=c]{90}{MRW 2D}}
 & \boldmath{$ 1.5$}    &   $0.84$   &   $0.68$   &   $0.51$  & $5 \,/\, 1$   &   $4 \,/\, 1$   &   $3 \,/\, 1$ \\
 & \boldmath{$ 2.8$}    &   $0.76$   &   $0.59$   &   $0.42$  & $5 \,/\, 1$   &   $4 \,/\, 1$   &   $3 \,/\, 1$ \\
 & \boldmath{$ 5.3$}    &   $0.74$   &   $0.51$   &   $0.22$  & $5 \,/\, 1$   &   $4 \,/\, 1$   &   $3 \,/\, 1$ \\
 & \boldmath{$ \infty$} &   $0.46$   &   $0.15$   &   $-0.17$ & $5 \,/\, 2$   &   $3 \,/\, 2$   &   $1 \,/\, 1$ \\
\hline
\end{tabular}
\end{table}

\begin{table}[h!]
\caption{\label{tab:cm}$\log_{10}(rmse)$ of $\cpfinico{m}$ and $\cpcoe{m}$ for MRW with different values of $p_0$
  and  $p$.
}
\begin{center}
\begin{tabular}{*{7}{c}}
\hline
& $p_0$ & $p=0.5$ &
 $p=2$ & $p=5$ & $p=\infty$ & DWT \\
\hline
\multirow{4}{*}{$c(1)$}
 & $ 1.3$    & $-1.76$   &      $-1.61$   &   $-1.48$   &   $-0.513$ & $-1.77$ \\
 & $ 2.5$    & $-1.8$   &      $-1.66$   &   $-1.52$   &   $-0.593$ & $-1.79$ \\
 & $   5$    & $-1.78$   &   $-1.65$   &   $-1.51$   &   $-0.807$ & $-1.68$ \\
 & $ \infty$ & $-1.79$   & $-1.78$   &   $-1.76$   &   $-1.75$ & $-1.69$ \\
\hline
\multirow{4}{*}{$c(2)$}
 & $ 1.3$    & $-2.16$   &    $-2.13$   &   $-2.04$   &   $-1.92$ & $-1.37$ \\
 & $ 2.5$    & $-2.16$   &    $-2.17$   &   $-2.07$   &   $-1.89$ & $-1.35$ \\
 & $   5$    & $-2.16$   &     $-2.1$   &   $-2$   &   $-1.9$ & $-1.33$ \\
 & $ \infty$ & $-2.05$   &    $-1.9$   &   $-1.84$   &   $-1.85$ & $-1.38$ \\
\hline
\multirow{4}{*}{$c(3)$}
 & $ 1.3$    &  $-2.1$   &      $-2.09$   &   $-2.02$   &   $-1.83$ & $-0.727$ \\
 & $ 2.5$    &  $-2.05$   &    $-1.96$   &   $-1.85$   &   $-1.74$ & $-0.685$ \\
 & $   5$    &  $-2.07$   &     $-2$   &   $-1.88$   &   $-1.78$ & $-0.688$ \\
 & $ \infty$ &  $-1.88$   &    $-1.7$   &   $-1.65$   &   $-1.63$ & $-0.625$ \\
\hline
\end{tabular}
\end{center}
\end{table}

\section{Finite-resolution effect in heart rate data}
\label{sec:Illustration}

Finally, the impacts of finite-resolution effects and the importance of using corrected $p$-leader
scaling exponent analysis is illustrated on heart rate (HR) analysis of Normal Sinus Rhythm, made
available by Physionet \cite{physionet}.
In this database, heart beats (RR intervals) were extracted by a standard automated procedure and
revised by experts. Following standard practice,  RR intervals were interpolated into a regularly
sampled time series, using cubic splines, at $f_s=4$ Hz.

\begin{figure}[t]
\centering
\includegraphics{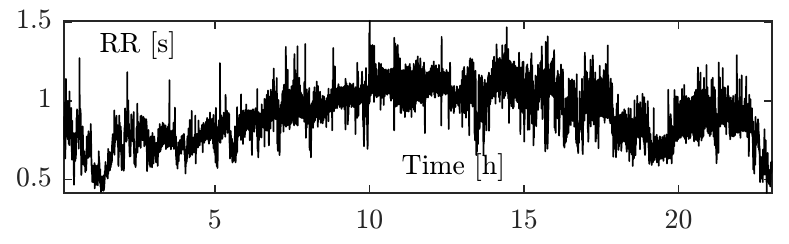}
\caption{\textbf{Sample heart rate data.} Record \emph{nsr046} of the  Normal Sinus Rhythm Physionet Database. }
\label{fig:hrv_data}
\end{figure}

\begin{figure}[t]
\centering
\begin{tabular}{cc}
\raisebox{-0.5\height}{\hspace*{-4mm}\includegraphics{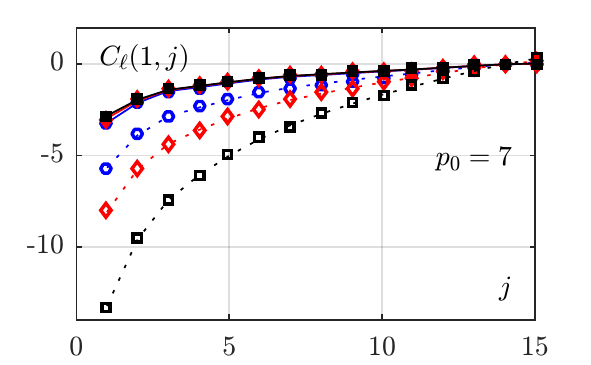}} &
\raisebox{-0.5\height}{\hspace*{-7mm}\includegraphics[scale=0.6]{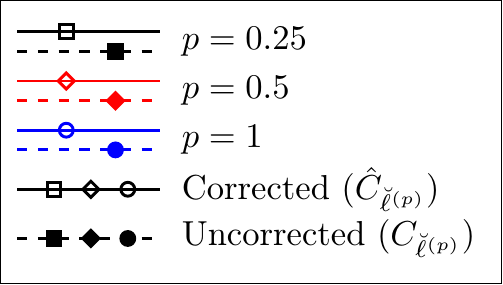}}
\end{tabular}
\caption{\textbf{Heart rate data: logscale diagrams.} Solid lines represent corrected
  $\cumjleadfinico{1}$, while dashed lines with solid markers represent uncorrected
  $\cumjleadfini{1}$. Marker-styles and colors indicate different values of $p$.
  The $\cumjleadfinico{1}$ coincide for different $p$,  as opposed to  uncorrected
  $\cumjleadfini{1}$.
}
\label{fig:hrv}
\end{figure}

The time series corresponding to record \emph{nsr046} is shown in Fig.~\ref{fig:hrv_data}.
The procedures described in Sec.~\ref{sec:min_reg_hyp} enable us to estimate $\hat{p}_0=7$.
Fig.~\ref{fig:hrv} reports the
corrected $\cumjleadfinico{1}$ and uncorrected $\cumjleadfini{1}$ cumulants for $p=0.25, 0.5$ and $1$.
Fig.~\ref{fig:hrv} clearly shows that uncorrected estimates for all chosen $p$s are affected by  finite-resolution effects.
Conversely, corrected estimates $\cumjleadfinico{1}$ for all $p$ collapse to a single function
$C(1,j)$ which can hence be considered as the actual scaling behavior of these data.
This example  illustrates that the estimation of the position of the location of the maximum
of the multifractal spectrum, which has been shown to be a relevant feature to discriminate healthy
from nonhealthy HR \cite{Spilka2016jbhi}, is  biased by finite-resolution effects, that can be well
accounted for by the proposed systematic correction.

\section{Conclusions}
\label{sec:Conclusions}

This contribution reports a thorough analysis of the finite-resolution effect that arises
when computing $p$-leaders from finite-resolution data.
Explicit closed-form relations were derived to account for such finite-resolution effect,
permitting to define corrected estimators than can be efficiently used in practice.
The complicated nonlinear definition of  $p$-leaders precluded their closed-form computation for
general functions or random processes; in consequence, no general proof of the validity of the
proposed corrections has been obtained so far ---but is undergoing further investigation.
Nonetheless, we assessed their effectiveness in two ways.
First, a theoretical analysis of multifractal cascades enabled an explicit computation of
$p$-leaders  and showed the relevance of the proposed corrections for both 1D and 2D cascades, with
very different correlation structures.
Second, numerical simulations allowed to show that the proposed correction is valid for several
types of multifractal processes of different natures.
Further, it was shown that corrected $p$-leaders have better estimation performance than wavelet
coefficients or state-of-the-art wavelet leaders.
Finally, the relevance of these issues for real-life heart rate data was illustrated.
The developments proposed in this work permit  to make use of the theoretical
and practical benefits of $p$-leaders for the multifractal analysis of data.
A {\sc Matlab} toolbox for $p$-leader multifractal analysis is
available at {\tt http://www.irit.fr/\%7EHerwig.Wendt/}.

\appendix
\label{app:proofs}

\subsection{Proof of Proposition \ref{prop:corr_cums}}
\label{app:corr_cums}

We follow the derivation of the log-cumulants in \cite{pLeadersPartII2015,Castaing1993,Wendt2007d}.
For infinite-resolution $p$-leaders, assuming that the moments of order $q$ exist and that
$\EE{(\rlead)^{q}}=F_q2^{-j\zetaq(q)}$, a standard generating function argument yields, for $q$
close to $0$, that
\begin{equation}
\label{eq:2}
F_q2^{-j\zetaq(q)}=\log\EE{e^{q\rlead}}=\sum_{m\geq 1}\cumjlead{m}\frac{q^m}{m!}.
\end{equation}

Now we consider finite-resolution $p$-leaders. Assuming  (\ref{eq:sf-nonlinear-corr}) and
(\ref{eq:sf-corr-equivalence}), and that moments of order $q$ exist,  we can deduce that the
expectation of $p$-leaders satisfies
\begin{equation}
\label{eq:exp_corr_lead}
\EE{(\rleadfini)^{q}}\corr{-\frac{q}{p}}{\eta(p)} = F_q 2^{-j\zetaq(q)}.
\end{equation}
Then, for $q$ close to $0$ we have
\begin{align}
\label{eq:4}
F_q 2^{-j\zetaq(q)}
&= \log\EE{e^{q\rleadfini}}-\frac{q}{p}\log\corr{}{\eta(p)}\\
&= \sum_{m\geq 1}\cumjleadfini{m}\frac{q^m}{m!} - \frac{q}{p}\log\corr{}{\eta(p)}\\
&= \sum_{m\geq 1}\cumjleadfinico{m}\frac{q^m}{m!},
\end{align}
where
\begin{equation}
\cumjleadfinico{m} =
\begin{cases}
\cumjleadfini{m} - \frac{q}{p}\log\corr{}{\eta(p)} &\text{if } m=1\\
\cumjleadfini{m} & \text{if } m\geq 2
\end{cases},
\end{equation}
which proves the direct statement. The  proof of the converse  statement is similar.

\subsection{Proof of (\ref{eq:eta_dbwc_2d})}  
\label{app:dbwc_eta}
Coefficients $d_{j,k_1,k_2}$ take
on values of the form $w_0^{n_0}w_1^{n_1}w_2^{n_2}w_3^{n_3}$, with $n_0+n_1+n_2+n_3=j$. Further,
from the tree structure of the cascade we have that
$$
\#\left\{ d_{j,\cdot,\cdot} \::\: n_0+n_1+n_2+n_3=j  \right\} =
\binom{j}{n_0,n_1,n_2,n_3}.
$$
Therefore, the wavelet structure function can be computed by application of  the multinomial theorem:
\begin{align}
\nonumber
\sfcoef &=2^{-2j}\sum_{k_1,k_2} \sum_i|\coe{j,k_1,k_2}^{(i)}|^q\\
&= \Vert\alpha\Vert_q^q \left( \frac{w^q_0+w^q_1+w^q_2+w^q_3}{4}\right)^j.
\end{align}
The wavelet scaling function $\eta$ is defined by the scaling relation $S_c(q,j) \sim
K2^{-j\eta(q)}$. Therefore, we can define:
\begin{equation}
\label{eq:eta2d}
\eta(q) = 2 -\log_2\sum_{m=0}^3w_m^q,\quad\text{for } q>0.
\end{equation}

\subsection{Proof of Proposition \ref{prop:dbwc_est}}
\label{app:dbwc_prop}
The  structure of the cascade implies that, for a fixed point $x$ and $j'>j$,
$d_{\lambda'(x)}=d_{\lambda(x)}w_{m_1}w_{m_2}\cdots w_{m_{j'-j}}$,
where the $m_i\in\{0,1,2,3\}$ take on values depending on each particular path on the subtree rooted in $d_{\lambda}$.
Then, using \eqref{eq:eta2d} we can compute the restricted $p$-leaders as

\begin{align}
\nonumber
\rleadfini &= \left( \sum_{\lambda'\subset\lambda}\sum_i|\coe{\lambda'}^{(i)}|^p2^{-2(j'-j)} \right)^{1/p}\\
\nonumber
&= \Vert\alpha\Vert_p d_{\lambda}\left(\sum_{l=0}^{\bar j- j}2^{-2l}\sum_{k_1,k_2} d_{j+l,k_1,k_2}^p\right)^{1/p}\\
\label{eq:pl-dwbc}
&= \Vert\alpha\Vert_p d_{\lambda} \left( \sum_{l=0}^{\bar j-j} 2^{-l\eta(p)} \right)^{1/p}.
\end{align}

Using (\ref{eq:def_sf}) and \eqref{eq:pl-dwbc}, we get:
\begin{equation}
\label{eq:sf-dwbc-pre}
\sfleadfini = \Vert\alpha\Vert_p \frac{1}{n_j}\sum_{k_1,k_2} d^q_{\lambda} \left( \sum_{l=0}^{\bar j-j} 2^{-l\eta(p)} \right)^{q/p}.
\end{equation}

For finite-resolution  $\bar j<\infty$,  the geometric sum in \eqref{eq:sf-dwbc-pre} adds up to
$\frac{1-2^{-(\bar j- j +1)\eta(p)}}{1-2^{-\eta(p)}}$ and hence the structure function is:
\begin{equation}
\label{eq:sf-dbwc-app}
\sfleadfini = \sfcoef \Vert\alpha\Vert_p  \left( \frac{1-2^{-(\bar j- j  +1)\eta(p)}}{1-2^{-\eta(p)}}\right)^{q/p}.
\end{equation}


\subsection{Proof of Proposition \ref{prop:rws-sf-nn}}  
\label{app:rws-sf-nn}

Let $n=q/p\in\NN$, and let $l=j'-j$ and $m=\bar j - j$.
The expected $p$-leader  structure function is $\EE{\sfleadfini}=\EE{\left( \rleadfini \right)^q}$. Further,
\begin{align}
\label{ec:rws-qth-moment}
\EE{\left(\rleadfini\right)^q}&=\EE{ \left( \sum_{l =0}^{m} \sum_{k=1}^{2^l}\coe{j+l,k}^p 2^{-l} \right)^{n} }\\
\label{eq:rws-stilde}
  &= \EE{ \left( \sum_{l=0}^m \tilde S_\lambda (p,l) \right)^{n} },
\end{align}
where $\tilde S_\lambda$ is the structure function of the RWS subtree rooted at coefficient $\coelam$.
Construction of the RWS implies that
\begin{equation}
\label{eq:ee-sf}
\EE{\tilde S_{\lambda}(p,l)}=\EE{\coe{j+l,\cdot}^p}=\EE{\{W^p\}^{(j+l)}}=2^{-(j+l)\eta(p)}.
\end{equation}

\paragraph*{Lower bound}
Since the function $x\mapsto x^n$ is convex for $n\in\NN^+$, we use Jensen's  inequality to get
\begin{align}
\label{eq:6}
\EE{ \left( \sum_{l=0}^m \tilde S_{\lambda} (p,l) \right)^{n} } &\geq
 \left( \sum_{l=0}^m \EE{ \tilde S_{\lambda} (p,l)} \right)^{n} \\
&\geq \left( \sum_{l=0}^m 2^{-(j+l)\eta(p)} \right)^n\\
&\geq 2^{-jn\eta(p)}\gamma^n(j, \eta(p)),
\end{align}
which proves the lower bound in (\ref{eq:rws-sf-nn}).

\paragraph*{Upper bound}
First we use the multinomial theorem on  (\ref{eq:rws-stilde}) and the fact that $\tilde
S_{\lambda}(p,l)$ are independent wrt $l$:
\begin{equation}
\label{eq:multin-first}
\EE{ \left( \sum_{l=0}^m \tilde S_{\lambda} (p,l) \right)^{n} } =
\sum_{\sum r_l=n}\binom{n}{r_0,\ldots,r_m} \prod_{l=0}^m \EE{\tilde S^{r_l}_{\lambda}(p,l)},\\
\end{equation}

Since  $r_l\in\NN$ for all $l$, we use (\ref{eq:ee-sf}) and the multinomial theorem again, which reads for $r\in\NN$:
\begin{equation}
\label{eq:multin-second}
\EE{ \tilde S^r_{\lambda} (p,l)} = 2^{-rl}\sum_{\sum s_k=r}\binom{r}{s_1,\dots,s_{2^l}} \prod_{k=1}^{2^l} 2^{-(j+l)\eta(p\,s_k)},
\end{equation}
where we have used the independence of $\coe{j+l,k}$ and \eqref{eq:ee-sf}.
Since $s_k\in\NN$ for all $k$, $\eta$ is concave and $\eta(0)=0$ we have:
\begin{equation}
\label{eq:bound-eta-upper}
\eta(p\,s_k)\geq \frac{s_{k}}{r}\eta(rp).
\end{equation}
Using (\ref{eq:bound-eta-upper}) in (\ref{eq:multin-second}) yields
\begin{equation}
\label{eq:1}
\EE{ \tilde S_{\lambda}^r (p,l)} \leq 2^{-(j+l)\eta(rp)}.
\end{equation}
Finally, (\ref{eq:1}) in (\ref{eq:multin-first}) produces
\begin{equation}
\label{eq:proof-rws-sf-nn}
\EE{ \left( \sum_{l=0}^m \tilde S_{\lambda} (p,l) \right)^{n} } \leq
\EE{\coelam^q}\gamma^n\left(j,\frac{\eta(np)}{n}\right),
\end{equation}
which proves the upper bound in (\ref{eq:rws-sf-nn}).

\subsection{Proof of Consequence \ref{prop:rws-sf-rr}}  
\label{app:rws-sf-rr}

Let $r\in\RR^+$ and $n\in\NN^+$.
Consider the polynomial interpolation
\begin{equation}
\label{eq:interp}
\sfleadfini[rp] = \sum_{n=0}^N \sfleadfini[np] h_n(x),
\end{equation}
where the $h_n(x)$ are the Lagrange basis polynomials.
Under the  smoothness assumption (i.e. $\Vert\partial^{N+1} \sflead[rp]/\partial
r^{N+1}\Vert_{\infty}$ is small enough) we can ignore the interpolation error.
Use of the lower bound in (\ref{eq:rws-sf-nn}) in (\ref{eq:interp}) yields
\begin{align}
\nonumber
\sfleadfini[rp] &\geq \sum_{n=0}^N \sfcoef[np] \,\corr{np}{\eta(p)}\,\corrl{n}{p} h_n(x),\\
\label{eq:proof-rws-rr-lower}
  &\geq \sfcoef[rp] \,\corr{rp}{\eta(p)}\,\corrl{r}{p},
\end{align}
which proves the lower bound. The  proof of the converse  statement is similar.

\subsection{Proof of Propositions \ref{prop:rwc_sf_p} and \ref{prop:rwc_sf_2p} }  
\label{app:rwc-sf-2p}

Proposition \ref{prop:rwc_sf_p} follows from  \eqref{ec:rws-qth-moment}, with
$n=1$. We  now prove Proposition \ref{prop:rwc_sf_2p}.
Expanding \eqref{ec:rws-qth-moment}, with $n=2$, the expected value of restricted $p$-leaders is given by
\begin{equation}
\label{eq:lead-moment-2}
\EE{(\rleadfini)^{2p}}=
\EE{\coelam ^{2p}}
\sum_{l_1=0}^{\bar j-j} \sum_{l_2=0}^{\bar j-j}\sum_{k_1=1}^{2^{l_1}}\sum_{k_2=1}^{2^{l_2}}
\EE{\coe{l_1,k_1}^p\coe{l_2,k_2}^p 2^{-l_1-l_2}}.
\end{equation}
Since wavelet coefficients of RWC are not independent, we cannot factorize the expected
value  $\EE{\coe{\lambda_1}\coe{\lambda_2}}$.

Let   $h:\ZZ^4\to [0,\min(l_1,l_2)]$ such that $h(l_1,l_2,k_1,k_2)$ is the  scale of the lowest
common ancestor between coefficients $\coe{l_1,k_1}$ and $\coe{l_2,k_2}$.
Then, as in \cite{Arneodo1998rwc}, we can write
$\coe{l_1,k_1} = W_1\ldots W_{h(l_1,l_2,k_1,k_2)} \, W_{h(l_1,l_2,k_1,k_2)+1}^{(1)}\ldots
W_{l_1}^{(1)}$ and
$\coe{l_2,k_2} = W_1\ldots W_{h(l_1,l_2,k_1,k_2)} \, W_{h(l_1,l_2,k_1,k_2)+1}^{(2)}\ldots W_{l_2}^{(2)}$
where $W_i$, $W^{(1)}_i$ and $W^{(2)}_i$ are iid random variables. Note
that the first ${h(l_1,l_2,k_1,k_2)}$
multipliers are the same for both coefficients. Therefore
\begin{equation*}
\EE{\coe{l_1,k_1}^p \coe{l_2,k_2}^p}=\EE{W^{2p}}^{h(l_1,l_2,k_1,k_2)}\EE{W^p}^{l_1+l_2-2h(l_1,l_2,k_1,k_2)}.
\end{equation*}

Using this in (\ref{eq:lead-moment-2}), and reordering we have:
\begin{multline}
\label{eq:18}
\EE{(\rleadfini)^{2p}} =
2^{-j\eta(2p)} \Biggl( \sum_{l_1 =0}^{\bar j-j}\sum_{l_2 =0}^{\bar j-j} 2^{-l_1-l_2}\EE{W^p}^{l_1+l_2}\cdot \\
\cdot\sum_{l=0}^{\min{(l_1,l_2)}} N(l_1,l_2,l)
\frac{\EE{W^{2p}}^{l}}{\EE{W^p}^{2l}} \Biggr).
\end{multline}
where the function $N$ represents the level sets of $h$:
\begin{equation}
\label{eq:rwc_def_N}
N(l_1, l_2, l) = \# \left\{ (k_1,k_2)\::\:h(l_1,l_2,k_1,k_2)=l \right\}.
\end{equation}

To compute $N$, let $l\leq\min(l_1,l_2)$. Consider the subtree rooted at coefficient $\coe{l,k}$ (note
that there are $2^l$ such subtrees): it  has $2^{l_2-l}$ children at scale $l_2$ and $2^{l_1-l}$ children
at scale $l_1$. Then, the total number of pairs $(\coe{l_1,k_1},\coe{l_2,k_2} )$ which have $\coe{l,k}$ as
a parent is $2^{l_1-l}2^{l_2-l}$.  Since there are $2^l$ possible choices for the root $\coe{l,k}$ we
have:
\begin{equation}
\label{eq:rwc_N}
N(l_1,l_2,l) = 2^{l_1-l}2^{l_2-l}2^l = 2^{l_1+l_2-l}
\end{equation}

Using (\ref{eq:rwc_N}) in (\ref{eq:18}), splitting the sum over $l$ and summing
the geometric sums, (\ref{eq:rwc-sf-2p}) follows.

\bibliographystyle{IEEEtran}
\bibliography{pLeaders3_Est_bib}
\end{document}